\documentclass{amsart}
\usepackage{amssymb}
\usepackage{latexsym}

\newcommand{\Q}{{\mathbb Q}}
\newcommand{\R}{{\mathbb R}}
\newcommand{\C}{{\mathbb C}}
\newcommand{\N}{{\mathbb N}}

\newcommand{\Z}{{\mathbb Z}}
\newcommand{\D}{{\mathbb D}}
\newcommand{\U}{{\mathbb U}}
\newcommand{\X}{{\mathbb S}}
\newcommand{\V}{{\mathbb L}}
\newcommand{\Y}{{\mathbb P}}
\newcommand{\G}{{\mathbb G}}

\newcommand{\CalW}{{\mathcal{W}}}
\newcommand{\CalA}{{\mathcal{A}}}

\newcommand{\Oomega}{(\Omega,T)}

\newtheorem{theorem}{Theorem}
\newtheorem{lemma}{Lemma}[section]
\newtheorem{prop}[lemma]{Proposition}

\begin{document}

\title{Uniform Szeg\H{o} Cocycles Over Strictly Ergodic Subshifts}

\author{David Damanik}

\address{Mathematics 253--37, California Institute of Technology,
Pasadena, CA 91125, U.S.A.}

\email{damanik@caltech.edu}

\thanks{D.\ D.\ was supported in part by NSF grant DMS--0500910}

\author{Daniel Lenz}

\address{Fakult\"at f\"ur Mathematik, TU Chemnitz, D-09107 Chemnitz, Germany}

\email{dlenz@mathematik.tu-chemnitz.de}

\begin{abstract}
We consider ergodic families of Verblunsky coefficients generated by minimal aperiodic
subshifts. Simon conjectured that the associated probability measures on the unit circle
have essential support of zero Lebesgue measure. We prove this for a large class of
subshifts, namely those satisfying Boshernitzan's condition. This is accomplished by
relating the essential support to uniform convergence properties of the corresponding
Szeg\H{o} cocycles.
\end{abstract}

\dedicatory{Dedicated to Barry Simon on the occasion of his 60th birthday.}

\maketitle

\section{Introduction}

$\Oomega$ is called a subshift over $\CalA$ if $\CalA$ is finite with discrete topology
and $\Omega$ is a closed $T$-invariant subset of $\CalA^\Z$, where $\CalA^\Z$ carries the
product topology and $T : \CalA^\Z \to \CalA^\Z$ is given by $(T s ) (n) = s (n+1)$.  A
function $F$ on $\Omega$ is called locally constant if there exists a non-negative
integer $N$ with
\begin{equation} \label{LocallyConstant}
F(\omega) = F(\omega') \quad \text{whenever} \quad   (\omega(-N),\ldots,
\omega(N))=(\omega'(-N),\ldots, \omega'(N)).
\end{equation}

The subshift $\Oomega$ is called minimal if every orbit $\{T^n \omega : n\in \Z\}$ is
dense in $\Omega$. It is called aperiodic if $T^n \omega\neq \omega$ for all $\omega\in
\Omega$ and $n\neq 0$.

Now, let a minimal subshift and   $f : \Omega \to \D$ continuous be given.  For $\omega
\in \Omega$, define $d\mu_\omega$ to be the probability measure on the unit circle
associated with the Verblunsky coefficients $\alpha_n(\omega) = f(T^n \omega)$, $n \ge
0$. Then, by minimality,  the essential support of the $d\mu_\omega$ does not depend on
$\omega$. It will be denoted by $\Sigma$.

In fact, $d\mu_\omega$ is the spectral measure of the CMV matrix $\mathcal C_\omega$
given by
$$
\begin{pmatrix}
{}& \bar\alpha_0(\omega) & \bar\alpha_1(\omega) \rho_0(\omega) & \rho_1(\omega)
\rho_0(\omega)
& 0 & 0 & \dots & {} \\
{}& \rho_0(\omega) & -\bar\alpha_1(\omega) \alpha_0(\omega) & -\rho_1(\omega)
\alpha_0(\omega)
& 0 & 0 & \dots & {} \\
{}& 0 & \bar\alpha_2(\omega) \rho_1(\omega) & -\bar\alpha_2(\omega) \alpha_1(\omega) &
\bar\alpha_3(\omega) \rho_2(\omega) & \rho_3(\omega) \rho_2(\omega) & \dots & {} \\
{}& 0 & \rho_2(\omega) \rho_1(\omega) & -\rho_2(\omega) \alpha_1(\omega) &
-\bar\alpha_3(\omega)
\alpha_2(\omega) & -\rho_3(\omega) \alpha_2(\omega) & \dots & {} \\
{}& 0 & 0 & 0 & \bar\alpha_4(\omega) \rho_3(\omega) & -\bar\alpha_4(\omega)
\alpha_3(\omega)
& \dots & {} \\
{}& \dots & \dots & \dots & \dots & \dots & \dots & {}
\end{pmatrix}
$$
Here, $\rho_n(\omega) = (1 - |\alpha_n(\omega)|)^{-1/2}$. Therefore, $\Sigma$ is the
essential spectrum of $\mathcal C_\omega$. In a completely analogous way, it is also
possible to define an extended CMV matrix, $\mathcal E_\omega$, acting on $\ell^2 (\Z)$
using the Verblunsky coefficients $\alpha_n(\omega) = f(T^n \omega)$, $n \in \Z$. Then,
$\Sigma$ is the spectrum of $\mathcal E_\omega$ for each $\omega\in \Omega$.

See Simon \cite{S,S2} for background on polynomials orthogonal with respect to a
probability measure on the unit circle (OPUC) and the associated Verblunsky coefficients
and CMV matrix. In \cite[Conjecture~12.8.2]{S2}, Simon conjectures the following:

\medskip

\noindent\textbf{Simon's Subshift Conjecture.} Suppose $\CalA$ is a subset of $\D$,  the
subshift $(\Omega,T)$ is minimal and aperiodic and let  $f : \Omega \to \D$, $f(\omega)=
\omega(0)$. Then, $\Sigma$ has zero Lebesgue measure.

\medskip

Our goal in this paper is to prove this conjecture (actually, in a stronger form) under a
mild assumption on the subshift that is satisfied in many cases of interest. We will give
some examples at the end of the paper and refer the reader to \cite{DL2} for a more
detailed discussion.

Before we formulate this assumption, let us recall that each subshift $\Oomega$ over
$\CalA$ gives rise to the associated set of words
$$
\CalW (\Omega) = \{ \omega(k) \cdots \omega(k + n -1) : k\in \Z, n\in \N, \omega \in
\Omega\}.
$$
For $w \in \CalW (\Omega)$, we define $V_w =\{\omega \in \Omega : \omega (1) \cdots
\omega (|w|) = w\}$.

$\Oomega$ is said to satisfy the Boshernitzan condition if it is minimal and there exists
an ergodic probability measure $\nu$ on $\Omega$, a constant $C>0$ and a sequence $(l_n)$
with $l_n \to \infty$ for $n \to \infty$ such that $|w| \nu(V_w)  \geq C $ whenever $w
\in \CalW (\Omega)$ with $|w|= l_n$ for some $n$; compare \cite{B3}.

\begin{theorem}\label{main}
Suppose the subshift $(\Omega,T)$ satisfies the Boshernitzan
condition. Let $f : \Omega \to \D$ be locally constant and
$\alpha_n(\omega) = f(T^n \omega)$ for $\omega \in \Omega$ and $n
\ge 0$. If for some $\omega \in \Omega$, $\alpha_n(\omega)$ is not
periodic, then for every $\omega \in \Omega$, the essential support
of the measure $d\mu_\omega$ associated with Verblunsky coefficients
$(\alpha_n(\omega))_{n \ge 0}$ has zero Lebesgue measure.
\end{theorem}

Note that by minimality of $\Omega$ and local constancy of $f$,
$\alpha_n(\omega)$ is periodic for some $\omega$ if and only if it
is periodic for every $\omega$.

Theorem~\ref{main} is the unit circle analogue of a result obtained in \cite{DL} for
discrete Schr\"odinger operators. Apart from addressing Simon's Subshift Conjecture, it
is our intention to introduce certain dynamical systems methods in the study of ergodic
Verblunsky coefficients and CMV matrices. These methods have been very successful in
recent studies of ergodic discrete Schr\"odinger operators, but have not yet found their
way into the OPUC literature to the degree they deserve.

\section{A Dynamical Characterization of the Essential Spectrum}

A continuous map $A : \Omega \to \G \V (2,\C)$ gives rise to a so-called cocycle, which
is a map from $\Omega \times \C^2$ to itself given by $(\omega,v) \mapsto (T
\omega,A(\omega)v)$. This map is usually denoted by the same symbol. When studying the
iterates of the cocycle, the following matrices describe the dynamics of the second
component:
$$
A(n,\omega) = \left\{\begin{array}{r@{\quad:\quad}l} A(T^{n-1} \omega) \cdots
A(\omega ) & n>0\\
Id & n=0\\ A (T^n \omega)^{-1} \cdots A (T^{-1}\omega)^{-1} & n < 0. \end{array} \right.
$$
By the multiplicative ergodic theorem, there exists a $\gamma(A)\in \R$ with
\begin{equation}\label{MultET}
\gamma (A)= \lim_{n \to \infty} \frac{1}{n} \log \| A(n,\omega)\|
\end{equation}
for $\mu$-almost every $\omega \in \Omega$. We say that $A$ is uniform if \eqref{MultET}
holds for all $\omega \in \Omega$ and the convergence is uniform on $\Omega$.

\begin{lemma}\label{L21}
Suppose $(\Omega,T)$ is a uniquely ergodic dynamical system. Let $A : \Omega \to \X \V
(2,\R)$ be continuous. Then the following are equivalent:
\begin{itemize}

\item[{\rm (i)}] $A$ is uniform and $\gamma (A)>0$.

\item[{\rm (ii)}] There exists a continuous map $P$ from $\Omega$ to the {\rm (}not
necessarily orthogonal{\rm )} projections on $\R^2$ and $C,\beta > 0$ such that
\begin{align*}
\| A(n,\omega) P(\omega) A(m,\omega)^{-1} \| & \le C e^{-\beta (n-m)} \quad n \ge m \\
\| A(n,\omega) (1 - P(\omega)) A(m,\omega)^{-1} \| & \le C e^{-\beta (m-n)} \quad n \le
m.
\end{align*}

\item[{\rm (iii)}] The exist continuous maps $U,V$ from $\Omega$ to the projective space
$\Y \R^2$ over $\R^2$  and $\tilde C, \tilde \beta > 0$ such that

\begin{align*}
\| A(n,\omega) u \| & \le \tilde C e^{-\tilde \beta n} \|u\| \quad \text{for}  \quad n \ge 0, \; u\in U(\omega)\\
\| A(-n,\omega) v \| & \le \tilde C e^{-\tilde \beta n} \|v\| \quad \text{for}  \quad n
\ge 0, \; v \in V(\omega).
\end{align*}

\end{itemize}
\end{lemma}
\begin{proof} The equivalence of (i) and (iii) is discussed in \cite[Theorem~4]{L2}.

The implication (ii) $\Longrightarrow$ (iii) follows from the case $m=0
$ in (ii) after one defines $U(\omega)= \mbox{Range}\, P(\omega)$ and
$V(\omega)= \mbox{Range} \, (1- P(\omega))$.

It remains to show (iii) $\Longrightarrow$ (ii). It is not hard to see that
$U(\omega)\neq V(\omega)$ for every $\omega \in \Omega$ (see \cite{L2} as
well). Thus, for $\omega\in \Omega$ fixed, any $x\in \R^2$ can be written
uniquely as
$$ x = u + v$$
with $u\in U(\omega)$ and $v\in V(\omega)$. We then define
$$ P(\omega) x := u.$$
Thus,
$$P(\omega) x\in U(\omega), \:\;\: (1- P(\omega)) x\in V(\omega)$$
for any
$\omega\in \Omega$ and $x\in \R^2$. As $U$ and $V$ are continuous, so is
$P$. Now, the case $m=0$ of (ii) follows directly. Moreover, as $U(T^m \omega)
= A(m,\omega) U(\omega)$ and $V (T^m\omega) = A(m,\omega)V(\omega)$, we have
$$ P(T^m \omega) A(m,\omega) = A(m,\omega) P(\omega)$$
and therefore
$$    A(m,\omega)^{-1} P(T^m \omega)  =  P(\omega) A(m,\omega)^{-1} $$
and similar with $P$ replaced by  $1 - P$.
Given this, the case of general $m$ in (ii) follows from the case $m=0$ with
$\omega$ replaced by $T^m \omega$.
\end{proof}

\noindent\textit{Remark.}  Condition (ii) is known as exponential
dichotomy. For a study of this condition in the context of Schr\"odinger operators
on the real line we refer the reader to \cite{J}.

\medskip

The cocycles we deal with do not take values in $\X\V(2,\R)$. They take values in $\X \U
(1,1)$, the set of $2\times 2$-matrices $A$ with determinant equal to one and $ A^{\ast}
J A = J$ for $J=\left(
\begin{array}{rr} 1 & 0 \\ 0 & -1 \end{array} \right)$.

\begin{lemma}\label{L22}
With the unitary matrix
$$
U = \frac{1}{\sqrt{2}} \left( \begin{array}{rr}  -i & 1 \\ i & 1 \end{array} \right)
$$
we have that $U^{-1} \, \X \U (1,1) \, U = \X \V (2,\R)$.
\end{lemma}

\begin{proof}
This follows from \cite[Proposition~10.4.1]{S2} and the discussion preceding it.
\end{proof}

\begin{theorem}\label{JL}
Suppose $(\Omega,T)$ is strictly ergodic and $f : \Omega \to \D$ is non-constant and
continuous. Define, for $z \in \partial \D$, the cocycle $A_z : \Omega \to \U (1,1)$ by
\begin{equation}\label{sc}
A_z(\omega) = (1 - |f(\omega)|^2)^{-1/2} \left( \begin{array}{cc} z & -
\overline{f(\omega)} \\ - f(\omega) z & 1 \end{array} \right).
\end{equation}
Then, $\partial \D \setminus \Sigma =   \{ z \in \partial \D : A_z \text{ is uniform and
} \gamma(A_z)
> 0 \}$.
\end{theorem}
\begin{proof} By Theorem~2.6 of \cite{GJ} (see also Remark~1 on page~143),
$$
\partial \D \setminus \Sigma =\{ z\in \partial \D : U^{-1} A_z U \text{ satisfies (ii)
of Lemma \ref{L21}}\}.
$$
By Lemma \ref{L21}, we can replace condition~(ii) by condition~(i) in the preceding
equality.  As $U$ is unitary, the assertion follows.
\end{proof}

\noindent\textit{Remark.} Since the matrices $A_z(n,\omega)$ are the transfer matrices
associated with the Szeg\H{o} recursion of the orthonormal polynomials with respect to
the measures $d\mu_\omega$, it is natural to call cocycles of the form \eqref{sc}
Szeg\H{o} cocycles.

\section{Boshernitzan's Condition and Uniform Convergence}

\begin{prop}\label{bosh}
Suppose the subshift $(\Omega,T)$ satisfies the Boshernitzan condition. Then, every
locally constant cocycle $A : \Omega \to \U (1,1)$ is uniform.
\end{prop}

\begin{proof}
Since $| \det A(\omega) | = 1$, it suffices to prove uniformity for the cocycle $\tilde A
: \Omega \to \X \U (1,1)$ given by $\tilde A(\omega) = (\det A(\omega))^{-1/2}
A(\omega)$. Here we choose an arbitrary but fixed branch of the square root. Note that
$\tilde A$ is locally constant since $A$ is locally constant. Obviously,
 $\tilde A$ is uniform if and only if the locally constant cocycle $\bar A = U^{-1} \tilde
A U$  with $U$ as in  Lemma~\ref{L22}  is uniform. Now, by Lemma~\ref{L22}, $\bar A$
takes values in $\X \L (2,\R)$.

It was shown in \cite[Theorem~1]{DL} that every locally constant $\X \V (2,\R)$ cocycle
over $(\Omega,T)$ is uniform if $(\Omega,T)$ satisfies the Boshernitzan condition.
\end{proof}

\section{Proof of the Main Result}

The OPUC version of Kotani theory yields the following consequence:
If the Verblunsky coefficients are ergodic, aperiodic, and take
finitely many values, then
\begin{equation}\label{kotani}
\mathrm{Leb} \left( \{ z \in \partial \D : \gamma(A_z)) = 0 \}
\right) = 0.
\end{equation}
See, for example, \cite[Theorem~10.11.4]{S2}. This result applies in the setting of
Theorem~\ref{main}.

\begin{proof}[Proof of Theorem~\ref{main}.]
By Theorem~\ref{JL}, the common essential support of the measures $d\mu_\omega$ is given
by
$$
\{ z \in \partial \D : \gamma(A_z) = 0 \} \cup \{ z \in \partial \D : \gamma(A_z) > 0
\text{ and $A_z$ is non-uniform} \}.
$$
The first set has zero Lebesgue measure by \eqref{kotani} and the second set is empty by
Proposition~\ref{bosh}.
\end{proof}

\section{Examples of Subshifts Satisfying Boshernitzan's Condition}

In this section we give a number of examples of aperiodic subshifts
that satisfy the Boshernitzan condition. As was shown in \cite{DL2},
most of the commonly studied minimal subshifts satisfy the
Boshernitzan condition. For example, subshifts obtained by codings
of rotations and interval exchange transformations.

Here we focus on subshifts associated with codings of rotations,
that is, models displaying a certain kind of quasi-periodicity. This
class has attracted a large amount of attention in other settings
(e.g., discrete Schr\"odinger operators and Jacobi matrices) but
their study in the context of OPUC is still in its early stages.

Let $\alpha \in (0,1)$ be irrational, $0 = \beta_0 < \beta_1 <
\cdots < \beta_{p-1} < \beta_p = 1$, and associate the intervals of
the induced partition with $p$ symbols $v_1,\ldots, v_p$:
$$
v_n(\theta) = v_k \Leftrightarrow \theta + n \alpha \, \mathrm{mod}
\, 1 \in [\beta_{k-1},\beta_k).
$$
We obtain a subshift over the alphabet $\{v_1,\ldots, v_p\}$,
$$
\Omega_{\alpha,\beta_1,\ldots,\beta_{p-1}} = \overline{ \{ v(\theta)
: \theta \in [0,1) \} }.
$$

The case $p=2$ is often of special interest. The following theorem
summarizes the results obtained in \cite{DL2} for this case.

\begin{theorem}
Let $\alpha$ be irrational. \\
{\rm (a)} If $\beta = m\alpha + n \, \mathrm{mod} \, 1$ with integers $m,n$,
then $\Omega_{\alpha,\beta}$ satisfies {\rm (B)}.\\
{\rm (b)} If $\alpha$ has bounded partial quotients, then
$\Omega_{\alpha,\beta}$ satisfies {\rm (B)} for every $\beta \in
(0,1)$.\\
{\rm (c)} If $\alpha$ has unbounded partial quotients, then
$\Omega_{\alpha,\beta}$ satisfies {\rm (B)} for Lebesgue almost
every $\beta \in (0,1)$.
\end{theorem}

Recall that $\alpha$ has bounded partial quotients if and only if the sequence $\{a_n\}$
in the continued fraction expansion of $\alpha$,
$$
\alpha = \cfrac{1}{a_1 + \cfrac{1}{a_2 + \cfrac{1}{a_3 + \cdots}}} ,
$$
is bounded.

It was also shown in \cite{DL2} that in the case where $\alpha$ has unbounded partial
quotients, there always exists a $\beta \in (0,1)$ such that $\Omega_{\alpha,\beta}$ does
not satisfy {\rm (B)}. These results show that the validity of (B) is well understood for
codings of rotations with respect to a partition of the circle into two intervals.

In the general case, we have the following theorem, which was
derived in \cite{DL2} from a result of Boshernitzan \cite{B2}.

\begin{theorem}
Let $\alpha \in (0,1)$ be irrational and suppose that $\beta_1,
\ldots, \beta_{p-1} \in \Q$. Then the subshift
$\Omega_{\alpha,\beta_1,\ldots,\beta_{p-1}}$ satisfies the
Boshernitzan condition {\rm (B)}.
\end{theorem}

\end{document}